\documentclass[preprint,review,10pt]{elsarticle}
\usepackage{amsfonts}
\usepackage{amsmath,amssymb}
\usepackage{graphicx}
\usepackage{setspace}
\usepackage[left=2cm,top=2cm,right=2cm]{geometry}

\newtheorem{lemma}{Lemma}

\newtheorem{remark}{Remark}

\newtheorem{theorem}{Theorem}
\numberwithin{equation}{section}
\journal{...}
\begin{document}
\title{$\alpha$$\beta$-Statistical Convergence of Modified $q$-Durrmeyer Operators}

\author[label*,label1,label2]{Vishnu Narayan Mishra}
\ead{vishnu\_narayanmishra@yahoo.co.in; vishnunarayanmishra@gmail.com}

\author[label1,label3]{Prashantkumar Patel}
\ead{prashant225@gmail.com}

\address[label1]{Department of Applied Mathematics $\&$ Humanities,
S. V. National Institute of Technology, Ichchhanath Mahadev Dumas Road, Surat-395 007 (Gujarat), India}
\address[label2]{L. 1627 Awadh Puri Colony Beniganj, Phase-III, Opposite - Industrial Training Institute (I.T.I.), Ayodhya Main Road, Faizabad, Uttar Pradesh 224 001, India}
\address[label3]{Department of Mathematics,
St. Xavier's College(Autonomous), Ahmedabad-380 009 (Gujarat), India}

\fntext[label*]{Corresponding author}
\begin{abstract}
In this work, we investigate weighted $\alpha$$\beta$-Statistical
approximation properties of $q$-Durrmeyer-Stancu operators. Also,
give some corrections in limit of $q$-Durrmeyer-Stancu operators
defined in \cite{mishra2013short} and discuss their convergence
properties.
\end{abstract}
\begin{keyword} Durrmeyer operators; Korovkin type theorems; Rate of the
weighted $\alpha$$\beta$-statistical convergent \end{keyword}

\maketitle
\section{Introduction}

The concept of statistical convergence has been defined by Fast
\cite{Fast1951sur} and studied by many other authors. It is well
know that every statistically convergent sequence is ordinary
convergent, but the converse is not true, examples and some
related work can be found in
\cite{dalmanog2010statistical,mishra2013statistical,orkcu2013approximation,Patelmishra20131,lin2014statistical,mursaleen2012weighted}.
The idea $\alpha\beta$-statistical convergence was
introduced by Akt\u{u}glu in \cite{Gadziev1974The}   as follows:\\
Let $\alpha(n)$ and $\beta(n)$ be two sequences positive number
which satisfy the following conditions
\begin{itemize}
\item[(i)] $\alpha$ and $\beta$ are both non-decreasing,
\item [(ii)] $\beta(n) \geq \alpha(n)$,
\item[(iii)]$\beta(n)-\alpha(n)\to \infty$ as $n\to \infty$
\end{itemize}
and let $\wedge$ denote the set of pairs $(\alpha,\beta)$
satisfying (i)-(iii). For each pair $(\alpha,\beta)\in \wedge$, $0
<\gamma\leq 1$ and $K \in \mathbb{N}$, we define
$\delta^{\alpha,\beta}(K,\gamma)$ in the following way
$$ \delta^{\alpha,\beta} (K,\gamma) = \lim_{n\to \infty}  \frac{|K\cap P_n^{\alpha,\beta}|}{\left(\beta(n)-\alpha(n)+1\right)^{\gamma}},$$
where $P_n^{\alpha,\beta}$ in the closed interval
$[\alpha(n),\beta(n)]$. A sequence $x = (x_k)$ is said to be
$\alpha\beta$-statistically convergent of order $\gamma$ to $\ell$
or $S_{\alpha\beta}^{\gamma}$-convergent, if
$$\delta^{\alpha,\beta}\left(\{ k: |x_k-\ell|\leq \epsilon\},\gamma\right)
= \lim_{n\to \infty} \frac{\big|\{k\in P_n^{\alpha,\beta} :
|x_k-\ell| \geq
\epsilon\}\big|}{\left(\beta(n)-\alpha(n)+1\right)^{\gamma}}=0.$$
The concept of weighted $\alpha\beta$-statistically convergent was
developed  by Karakaya and Karaisa \cite{karakaya2015korovkin}.
Let $s = (s_k)$ be a sequence of non-negative real numbers such
that $s_0
> 0$ and $$S_n = \sum_{k\in P_n^{\alpha,\beta}}s_k\to \infty,
\textrm{ as  } n \to \infty \textrm{ and } z_n^\gamma(x)
=\frac{1}{S_n^{\gamma}}\sum_{k\in P_n^{\alpha,\beta}}s_kx_k.$$ A
sequence $x = (x_k)$ is said to be weighted
$\alpha\beta$-statistically convergent of order $\gamma$ to $\ell$
or $S_{\alpha\beta}^{\gamma}$-convergent, if for every
$\epsilon>0$
$$\delta^{\alpha,\beta} \left(\{ k : s_k |x_k-\ell|\geq \epsilon\},\gamma\right)=\lim_{n\to \infty} \frac{1}{S_n^{\gamma}}|\{k\leq S_n : s_k |x_k-\ell| \geq \epsilon\}| =0$$
and denote $st_{\alpha\beta}^{\gamma} -\lim x = \ell$ or $x_k \to
\ell[\bar{S}_{\alpha\beta}^{\gamma}]$, where
$\bar{S}_{\alpha\beta}^{\gamma}$ denotes the set of all weighted
$\alpha\beta$-statistically convergent sequences of order
$\gamma$. \\
\indent The $q$-Bernstein operators were introduced by Phillips in
\cite{phillips1996bernstein} and they generalized the well-known
Bernstein operators. A survey of the obtained results and
references concerning $q$-Bernstein operators can be found in
\cite{ostrovska2007first}. It is worth mentioning that the first
generalization of the Bernstein operators based on $q$-integers
was obtained by Lupa\c{s} \cite{lupas1987q}. The Durrmeyer type
modification of $q$-Bernstein operators were established by Gupta
\cite{gupta2008rate} and it's local approximation, global
approximation and simultaneous approximation properties were
discussed in \cite{finta2009approximation}, we refer some of the
important papers in this direction as
\cite{mishra2014durrmeyer,buyukyazici2010approximation,vedi2015chlodowsky,mursaleen2015approximation,mursaleen2014statistical,edely2013approximation,mursaleen2012weighted,acar2014power,acar2015asymptotic}.
Also, better approximation properties were established by Gupta
and Sharma \cite{gupta2011recurrence}. Stancu type generalization
of the $q$-Durrmeyer operators were discussed by Mishra and Patel
\cite{mishra2013short,mishra2014generalized}, which define for
$f\in C\left([0,1]\right)$ as
\begin{eqnarray}\label{18.1.eq1}
D_{n,q}^{\varpi,\vartheta} &=& [n+1]_q \sum_{k=0}^{\infty} q^{-k}
p_{nk}(q;x) \int_0^1 f\left(\frac{[n]_qt+\varpi}{[n]_q+
\vartheta}\right)p_{nk}(q,qt) d_qt\\
&=& \sum_{k=0}^{\infty} A_{n,k}^{\varpi,\vartheta}(f)p_{nk}(q;x);
0\leq x\leq 1,\nonumber
\end{eqnarray}
where $\displaystyle p_{nk}(q;x) = {\binom{n}{k}}_q x^k
(1-x)_q^{n-k}$. We have used notations of $q$-calculus as given in
\cite{kac2002quantum}. Along the paper, $C ([a,b])$ denote by set
of continuous functions on interval $[a,b]$ and
$\|h\|_{C([a,b])}$ represents the sup-norm of the function
$h\mid_{[a,b]}$.\\
 \indent In this work, we establish
$\alpha\beta$-statistical convergence for operators
\eqref{18.1.eq1}. Also, in section \ref{18.sec3}, we discuss
convergence results of limit of $q$-Durrmeyer-Stancu operators
\eqref{18.1.eq1}.
\begin{lemma}[\cite{mishra2013short}]\label{l1}
We have\\
$D_{n,q}^{\varpi,\beta}(1; x)=1, \text{  } D_{n,q}^{\varpi,\vartheta}(t; x)= \displaystyle \frac{[n]_q+\varpi[n+2]_q+qx[n]_q^2}{[n+2]_q\left([n]_q+\vartheta\right)}$\\
and\\ $ D_{n,q}^{\varpi,\vartheta}(t^2; x)= \displaystyle
\frac{q^3[n]_q^3\left([n]_q-1\right)x^2+\left(\left(q(1+q)^2+2\varpi
q^4\right)[n]_q^3+2\varpi
q[3]_q[n]_q^2\right)x}{([n]_q+\vartheta)^2[n+2]_q[n+3]_q}+
\frac{(1+q+2\varpi
q^3)[n]_q^2+2\varpi[3]_q[n]_q}{([n]_q+\vartheta)^2[n+2]_q[n+3]_q}\\~~~~~~~~~~~~~~~~~~~~~+\frac{\varpi^2}{([n]_q+\vartheta)^2}.$
\end{lemma}
\begin{remark}\label{l4}
By simple computation, we can find the central moments \\
$\displaystyle \delta_{n}(x)=D_{n, q}^{\varpi, \vartheta}(t-x;x)= \left(\frac{q[n]_q^2}{[n+2]_q([n]_q+\vartheta)}-1\right)x+\frac{[n]_q+\varpi[n+2]_q}{[n+2]_q([n]_q+\vartheta)},$\\
$\displaystyle \gamma_n(x) =
 D_{n, q}^{\varpi, \vartheta}((t-x)^2;x)= \frac{q^4[n]_q^4-q^3[n]_q^3-2q[n]_q^2[n+3]_q([n]_q+\vartheta)+[n+2]_q[n+3]_q([n]_q+\vartheta)^2}{([n]_q+\vartheta)^2[n+2]_q[n+3]_q}x^2\\
\indent \indent \indent \indent \indent \indent \indent +\frac{q(1+q)^2[n]_q^3+2q\varpi [n]_q^2[n+3]_q-\left(2[n]_q+2\varpi[n+2]_q\right)[n+3]_q([n]_q+\vartheta)}{([n]_q+\vartheta)^2[n+2]_q[n+3]_q}x\\
\indent \indent \indent \indent \indent \indent \indent
+\frac{(1+q)[n]_q^2+2\varpi[n]_q[n+3]_q}{([n]_q+\vartheta)^2[n+2]_q[n+3]_q}.$
\end{remark}
\section{$\alpha\beta$-Statistical Convergence}
\begin{theorem}[\cite{karakaya2015korovkin}]\label{18.2.thm2} Let $(L_k)$ be a sequence of positive linear
operator from $C\left([a, b]\right)$ into $C\left([a, b]\right)$.
Then for all $f \in C\left([a, b]\right)$
$$\bar{S}_{\alpha\beta}^{\gamma} - \lim_{k\to \infty} \|L_k(f, x)-f(x) \|_{C\left([a, b]\right)}=0$$
if and only if
$$\bar{S}_{\alpha\beta}^{\gamma} - \lim_{k\to \infty} \|L_k(x^i, x)-x^i \|_{C\left([a, b]\right)}=0,~~~i=0,1,2.$$
\end{theorem}
Let $\{q_n\}$ be a sequence in the interval $[0,1]$ satisfying
\begin{equation}\label{18.2.eq1}
\bar{S}_{\alpha\beta}^{\gamma} - \lim_{k\to \infty}
q_n=1,~~ \bar{S}_{\alpha\beta}^{\gamma} - \lim_{k\to \infty}
(q_n)^n=a\in (0,1),~~\bar{S}_{\alpha\beta}^{\gamma} - \lim_{k\to
\infty} \frac{1}{[n]_q}=1
\end{equation}

\begin{theorem}\label{18.2.thm3}
Let $\{q_n\}$ be a sequence satisfying (\ref{18.2.eq1}) and
$D_{n,q}^{\vartheta,\varpi}$ as defined in \eqref{18.1.eq1}. For
any $f\in C\left([0,1]\right)$, we have
$$\bar{S}_{\alpha\beta}^{\gamma} -\lim_{k\to \infty} \|D_{n,q}^{\vartheta,\varpi}(f, x)-f(x) \|_{C\left([0, 1]\right)}=0.$$
\end{theorem}
\textbf{Proof: } By Theorem \ref{18.2.thm2}, it is enough to prove
that
\begin{equation}
\bar{S}_{\alpha\beta}^{\gamma} -\lim_{k\to \infty}
\|D_{n,q}^{\vartheta,\varpi}(t^j,
x)-x^j\|_{C\left([0,1]\right)}=0,~~~j=0,1,2
\end{equation}
From the $D_{n,q}^{\vartheta,\varpi}(1, x)=1$, it is easy to
obtain that
$$\bar{S}_{\alpha\beta}^{\gamma} -\lim_{k\to \infty} \|D_{n,q}^{\vartheta,\varpi}(1, x)-1 \|_{C\left([0, 1]\right)}=0.$$
Now,
\begin{eqnarray*}
|D_{n,q}^{\varpi,\vartheta}(t;
x)-x|&\leq &\bigg|\frac{q[n]_q^2-[n+2]_q([n]_q+\vartheta)}{[n+2]_q([n]_q+\vartheta)}\bigg|+\bigg|\frac{[n]_q+\varpi[n+2]_q}{[n+2]_q([n]_q+\vartheta)}\bigg|\\
&=&\bigg|\frac{[n]_q\left(q[n]_q-[n+2]_q\right)-\vartheta[n+2]_q}{[n+2]_q([n]_q+\vartheta)}
\bigg|+\bigg|\frac{[n]_q+\varpi[n+2]_q}{[n+2]_q([n]_q+\vartheta)}\bigg|\\
&\leq
&\bigg|\frac{[n]_q(1+q^{n+1})}{[n+2]_q([n]_q+\vartheta)}\bigg|+\bigg|\frac{\vartheta}{[n]_q+\vartheta}
\bigg|+\bigg|\frac{[n]_q+\varpi[n+2]_q}{[n+2]_q([n]_q+\vartheta)}\bigg|\\
\end{eqnarray*}
Using equation (\ref{18.2.eq1}), we get
$$\bar{S}_{\alpha\beta}^{\gamma} -\lim_{k\to \infty} \frac{[n]_q(1+q^{n+1})}{[n+2]_q([n]_q+\vartheta)}=0;~~~
 \bar{S}_{\alpha\beta}^{\gamma} -\lim_{k\to \infty} \bigg|\frac{\vartheta}{[n]_q+\vartheta}
\bigg|=0 $$ and
$$\bar{S}_{\alpha\beta}^{\gamma} -\lim_{k\to
\infty}\frac{[n]_q+\varpi[n+2]_q}{[n+2]_q([n]_q+\vartheta)}=0$$
Define the following sets:
$$A=\left\{n\in \mathbb{N} :\|D_{n,q}^{\varpi,\vartheta}(\cdot;
x)-x\|_{C\left([a, b]\right)}\geq
\epsilon\right\};~~~A_1=\left\{n\in \mathbb{N}
:\|\frac{[n]_q(1+q^{n+1})}{[n+2]_q([n]_q+\vartheta)}\geq
\frac{\epsilon}{3}\right\};$$
$$A_2=\left\{n\in \mathbb{N}
:\|\frac{\vartheta}{[n]_q+\vartheta}\geq
\frac{\epsilon}{3}\right\},~~~~A_3=\left\{n\in \mathbb{N}
:\|\frac{[n]_q+\varpi[n+2]_q}{[n+2]_q([n]_q+\vartheta)}\geq
\frac{\epsilon}{3}\right\},$$ Then, we obtain $A\subset A_1 \cup
A_2\cup A_3$, which implies that
$\delta_{\gamma}^{\alpha,\beta}(A) \leq
\delta_{\gamma}^{\alpha,\beta}(A_1)+\delta_{\gamma}^{\alpha,\beta}(A_2)+\delta_{\gamma}^{\alpha,\beta}(A_3)$
and hence
$$\bar{S}_{\alpha\beta}^{\gamma} -\lim_{k\to \infty}
\|D_{n,q}^{\vartheta,\varpi}(t, x)-x\|_{C\left([0, 1]\right)}=0.$$
Similarly, we have
\begin{eqnarray*}
|D_{n,q}^{\varpi,\vartheta}(t^2; x)-x^2| &\leq&\bigg|
\frac{q^3[n]_q^3\left([n]_q-1\right)}{([n]_q+\vartheta)^2[n+2]_q[n+3]_q}-1\bigg|\\
&&+\bigg|\frac{\left(\left(q(1+q)^2+2\varpi
q^4\right)[n]_q^3+2\varpi
q[3]_q[n]_q^2\right)}{([n]_q+\vartheta)^2[n+2]_q[n+3]_q}\bigg|\\
&&+\bigg| \frac{(1+q+2\varpi
q^3)[n]_q^2+2\varpi[3]_q[n]_q}{([n]_q+\vartheta)^2[n+2]_q[n+3]_q}\bigg|+\bigg|\frac{\varpi^2}{([n]_q+\vartheta)^2}\bigg|\\
&\leq& \bigg|
\frac{q^3_n[n]_{q_n}^4(1-q_n^2)}{([n]_q+\vartheta)^2[n+2]_q[n+3]_q}\bigg|
+\bigg|\frac{\left(q(1+q)^2+2\varpi q^4\right)[n]_q^3}{([n]_q+\vartheta)^2[n+2]_q[n+3]_q}\bigg|\\
&&+\bigg|\frac{2\varpi
q[3]_q[n]_q^2}{([n]_q+\vartheta)^2[n+2]_q[n+3]_q}\bigg|+\bigg|
\frac{(1+q+2q^3\varpi)[n]_q^2}{([n]_q+\vartheta)^2[n+2]_q[n+3]_q}\bigg|\\
&&+\bigg|\frac{2\varpi
[3]_q[n]_q}{([n]_q+\vartheta)^2[n+2]_q[n+3]_q}\bigg|+\bigg|\frac{\varpi^2}{([n]_q+\vartheta)^2}
\bigg|.
\end{eqnarray*}
Again, using $\displaystyle \bar{S}_{\alpha\beta}^{\gamma} -
\lim_{k\to \infty} q_n=1,~~ \bar{S}_{\alpha\beta}^{\gamma} -
\lim_{k\to \infty} (q_n)^n=a\in
(0,1),~~\bar{S}_{\alpha\beta}^{\gamma} - \lim_{k\to \infty}
\frac{1}{[n]_q}=1$, we get
\begin{eqnarray*}
\bar{S}_{\alpha\beta}^{\gamma}-\lim_{n\to
\infty}\frac{q^3_n[n]_{q_n}^4(1-q_n^2)}{([n]_q+\vartheta)^2[n+2]_q[n+3]_q}&=&0,\\
\bar{S}_{\alpha\beta}^{\gamma}-\lim_{n\to \infty}\frac{\left(q(1+q)^2+2\varpi q^4\right)[n]_q^3}{([n]_q+\vartheta)^2[n+2]_q[n+3]_q}&=&0,\\
\bar{S}_{\alpha\beta}^{\gamma}-\lim_{n\to \infty}\frac{2\varpi
q[3]_q[n]_q^2}{([n]_q+\vartheta)^2[n+2]_q[n+3]_q}&=&0,\\
\bar{S}_{\alpha\beta}^{\gamma}-\lim_{n\to \infty}\frac{(1+q+2q^3\varpi)[n]_q^2}{([n]_q+\vartheta)^2[n+2]_q[n+3]_q}&=&0,\\
\bar{S}_{\alpha\beta}^{\gamma}-\lim_{n\to \infty}\frac{2\varpi
[3]_q[n]_q}{([n]_q+\vartheta)^2[n+2]_q[n+3]_q}&=&0,\\
\bar{S}_{\alpha\beta}^{\gamma}-\lim_{n\to
\infty}\frac{\varpi^2}{([n]_q+\vartheta)^2}&=&0.
\end{eqnarray*}
Now, consider the following sets:
\begin{eqnarray*}
B_1&:=&\left\{n\in \mathbb{N}: \frac{q^3_n[n]_{q_n}^4(1-q_n^2)}{([n]_q+\vartheta)^2[n+2]_q[n+3]_q}\geq \frac{\epsilon}{6}\right\},\\
B_2&:=&\left\{n\in \mathbb{N}:\frac{\left(q(1+q)^2+2\varpi q^4\right)[n]_q^3}{([n]_q+\vartheta)^2[n+2]_q[n+3]_q}\geq \frac{\epsilon}{6}\right\},\\
B_3&:=&\left\{n\in \mathbb{N}:\frac{2\varpi
q[3]_q[n]_q^2}{([n]_q+\vartheta)^2[n+2]_q[n+3]_q}\geq \frac{\epsilon}{6}\right\},\\
B_4&:=&\left\{n\in \mathbb{N}:\frac{(1+q+2q^3\varpi)[n]_q^2}{([n]_q+\vartheta)^2[n+2]_q[n+3]_q}\geq \frac{\epsilon}{6}\right\},\\
B_5&:=&\left\{n\in \mathbb{N}:\frac{2\varpi
[3]_q[n]_q}{([n]_q+\vartheta)^2[n+2]_q[n+3]_q}\geq \frac{\epsilon}{6}\right\},\\
B_6&:=&\left\{n\in
\mathbb{N}:\frac{\varpi^2}{([n]_q+\vartheta)^2}\geq
\frac{\epsilon}{6}\right\}.
\end{eqnarray*}
Consequently, we obtain $B\subset B_1 \cup B_2 \cup B_3 \cup B_4
\cup B_5 \cup B_6$, which implies that $\displaystyle \delta(B)
\leq \sum_{i=1}^6 \delta(B_i)$. Hence, we get
$$\bar{S}_{\alpha\beta}^{\gamma} -\lim_{k\to \infty}
\|D_{n,q}^{\vartheta,\varpi}(t^2,
x)-x^2\|_{C\left([0,1]\right)}=0.$$ This completes the proof of
Theorem \ref{18.2.thm3}.
\section{Limit
$q$-Durrmeyer-Stancu operators}\label{18.sec3} The authors found
mistake in the proof part of \cite[Theorem
2]{mishra2014generalized}. In \cite[Sec.
4]{mishra2014generalized}, authors defined the operators
$D_{\infty,q}^{\vartheta,\varpi}$  \cite[Eq.
(4.2)]{mishra2014generalized}, which depend on $[n]_q$ was
mistaken.
 So, follow by \cite[Theorem 2]{mishra2014generalized} the proof part have some errors.
With this note we correctly define the operators and prove Theorem
2 of \cite{mishra2014generalized}.\\ 
Here, we define the limit $q$-Durrmeyer-Stancu operators \eqref{18.1.eq1} as:\\
 Let $q\in (0,1) $ be fixed and $x\in [0,1]$, the operators $D_{\infty,q}^{\vartheta,\varpi}(f;x)$ is defined by
 \begin{eqnarray}\label{e7}
  D_{\infty,q}^{\vartheta,\varpi}(f;x)&=& \frac{1}{1-q}\sum_{k=0}^{\infty} p_{\infty k}(q;x)q^{-k} \int_0^1 f\left(\frac{ t+(1-q)\vartheta}{1+(1-q)\varpi}\right) p_{\infty k}
  (q;qt)d_qt\\
 &=& \sum_{k=0}^{\infty} A_{\infty k}^{\vartheta,\varpi}(f)  p_{\infty k }(q;x)\nonumber.
  \end{eqnarray}
  Using the fact that (see \cite{il2002convergence}), we have
\begin{equation}
\sum_{k=0}^{\infty} p_{\infty k}(q;x)=1 , ~~\sum_{k=0}^{\infty} (1-q^k)p_{\infty k}(q;x)=x,
\end{equation} and
\begin{equation}
 \sum_{k=0}^{\infty} (1-q^k)^2p_{\infty k}(q;x) =x^2+(1-q)x(1-x).
 \end{equation}
  Using (\ref{18.1.eq1}) and (\ref{e7}), it is easy to prove that\\
$\displaystyle D_{\infty,q}^{\vartheta,\varpi}(1;x)=1 , \text{  }
D_{\infty,q}^{\vartheta,\varpi}(t;x)=  \frac{1+
q(x-1)+\vartheta(1-q)}{1+\varpi(1-q)},\\~~
 \displaystyle D_{\infty,q}^{\vartheta,\varpi}(t^2;x)=\frac{ q^4 x^2+\left(q(1+q)(1-q^2)+2(1-q)q\vartheta\right)x}{(1+\varpi(1-q))^2}
 +\frac{((1+q)+2\vartheta+\vartheta^2)(1-q)^2}{(1+\varpi(1-q))^2}.
$\\
 For $f\in C[0,1], t>0$, we define the modulus of continuity $\omega(f,t)$ as follows:
 $$ \omega(f,t) = \sup \{ |f(x)-f(y)|:|x-y| \leq t, ~~x,y\in
 [0,1]\}.$$
 \begin{theorem}\label{t3}
 Let $0<q<1$ then for each $f\in C[0,1]$ the sequence $\{D_{n,q}^{\vartheta,\varpi}(f;x)\}$
 converges to $ D_{\infty,q}^{\vartheta,\varpi}(f;x)$ uniformly on $[0,1]$. Furthermore,
 $$ \displaystyle \|D_{n,q}^{\vartheta,\varpi}(f)-D_{\infty,q}^{\vartheta,\varpi}(f)\| \leq
 C_q^{\vartheta,\varpi}\omega(f,q^n).$$
  \end{theorem}
\textbf{Proof: } $  D_{\infty,q}^{\vartheta,\varpi}(f;x)$ and
$D_{n,q}^{\vartheta,\varpi}(f;x)$ reproduce constant function that
is $\displaystyle D_{n,q}^{\vartheta,\varpi}(1;x)=
D_{\infty,q}^{\vartheta,\varpi}(1;x)=1.$ Hence for all $x\in
[0,1)$, by definition of $D_{n,q}^{\vartheta,\varpi}(f;x)$ and
$D_{\infty,q}^{\vartheta,\varpi}(f;x)$, we know that
\begin{eqnarray*}
 |D_{n,q}^{\vartheta,\varpi}(f;x)-D_{\infty,q}^{\vartheta,\varpi}(f;x)|
 &=&
\bigg| \sum_{k=0}^{n} A_{n k}^{\vartheta,\varpi} (f) p_{nk}(q;x)-
\sum_{k=0}^{\infty} A_{\infty k}^{\vartheta,\varpi} (f) p_{\infty
k}(q;x)  \bigg|\\
 &=&  \bigg| \sum_{k=0}^{n} A_{n k}^{\vartheta,\varpi} \left(f-f(1)\right) p_{nk}(q;x)\\
 &&- \sum_{k=0}^{\infty} A_{\infty k}^{\vartheta,\varpi} \left(f-f(1)\right) p_{\infty k}(q;x)  \bigg|\\
 &\leq& \sum_{k=0}^{n} \big|A_{n k}^{\vartheta,\varpi} \left(f-f(1)\right)-A_{\infty k}^{\vartheta,\varpi} \left(f-f(1)\right)\big| p_{nk}(q;x)\\
&&+ \sum_{k=0}^{n} |A_{\infty k}^{\vartheta,\varpi}
\left(f-f(1)\right)| |p_{n k}(q;x)- p_{\infty k}(q;x)|\\
 &&+ \sum_{k=n+1}^{\infty} |A_{\infty k}^{\vartheta,\varpi} \left(f-f(1)\right)|p_{\infty k}(q;x)
 = I_1+I_2+I_3.
 \end{eqnarray*}
 By the well known property of modulus of continuity (see \cite{lorentz1953bernstein}), $ \omega(f,\lambda t) \leq (1+\lambda)\omega(f,t), \lambda >0,$
 we get\\
$$|f(t)-f(1)| \leq \omega(f, 1-t) \leq
\omega(f,q^n)\left(1+\frac{1-t}{q^n}\right).$$
 Thus
 \begin{eqnarray*}
 |A_{n k}^{\vartheta, \varpi}(f-f(1))| &=& \bigg| [n+1]_q \int_0^1 q^{-k} \left(f\left(\frac{[n]_qt+\vartheta}{[n]_q+\varpi}\right)-f(1)\right)p_{n k} (q;qt) d_qt\bigg|\\
 &\leq&  [n+1]_q \int_0^1 q^{-k} \bigg|f\left(\frac{[n]_qt+\vartheta}{[n]_q+\varpi}\right)-f(1)\bigg|p_{n k} (q;qt) d_qt\\
  &\leq&  [n+1]_q \int_0^1 q^{-k} \omega(f,q^n)\left(1+\frac{1}{q^n}\left(1-\frac{[n]_qt+\vartheta}{[n]_q+\varpi}\right)\right)p_{n k} (q;qt) d_qt\\
 &\leq& \omega(f,q^n)\left(1+ q^{-n}\left(1- \frac{[n]_q[k+1]_q-\vartheta[n+2]_q}{[n+2]_q([n]_q+\varpi)}\right)\right)
\end{eqnarray*}
 \begin{eqnarray*}
 |A_{n k}^{\vartheta, \varpi}(f-f(1))|   &\leq & \omega(f,q^n)\left(1+ \frac{q^{-n}[n]_q}{[n]_q+\varpi}\left(1-\frac{[k+1]_q}{[n+2]_q}\right)+ \frac{q^{-n}(\varpi-\vartheta)}{[n]_q+\varpi}\right)\\
 &=& \omega(f,q^n)\left(1+ q^{k+1-n}+ \frac{q^{-n}(\varpi-\vartheta)}{[n]_q+\varpi}\right).
 \end{eqnarray*}
Similarly,
 \begin{eqnarray*}
|A_{\infty k}^{\vartheta, \varpi}(f-f(1))|&=& \frac{q^{-k}}{1-q} \bigg| \int_0^1  \left(f\left(\frac{t+\vartheta(1-q)}{1+\varpi(1-q)}\right)-f(1)\right)p_{\infty k} (q;qt) d_qt \bigg|\\
 &\leq & \frac{q^{-k}}{1-q} \int_0^1 \omega(f,q^n)\left(1+\frac{1}{q^n}\left(1-\frac{t+\vartheta(1-q)}{1+\varpi(1-q)}\right)\right)p_{\infty k} (q;qt) d_qt\\
 &\leq&\frac{q^{-k}}{1-q} \int_0^1  \omega(f,q^n)\left(1+\frac{1}{q^n}(1-t)+\frac{1}{q^n}\frac{\varpi-\vartheta}{1+\varpi(1-q)}\right)p_{\infty k} (q;qt) d_qt\\
 &\leq& \omega(f,q^n)\left(1+ q^{k+1-n}+\frac{q^{-n}(\varpi-\vartheta)}{1+\varpi(1-q)}\right).
 \end{eqnarray*}
 From \cite[Eq.4.5]{mishra2014generalized}, we have
 \begin{equation}\label{e9*}
 |p_{n k}(q;x)-p_{\infty k}(q;x)|   \leq \frac{q^{n-k}}{1-q}\left( p_{n k}(q;x)+p_{ \infty k}(q;x)\right).
 \end{equation}
 Hence by using (\ref{e9*}), we have\\
\begin{eqnarray*}
&&| A_{n k}^{\vartheta, \varpi}(f-f(1))- A_{\infty k}^{\vartheta, \varpi}(f-f(1))| \\
 &\leq& [n+1]_q\int_0^1 q^{-k} \bigg|f\left(\frac{[n]_qt+\vartheta}{[n]_q+\varpi}\right)-f(1)\bigg|  p_{n k}(q;qt) d_qt + \frac{1}{1-q}\int_0^1 q^{-k} \bigg| f\left(\frac{t+\vartheta(1-q)}{1+\varpi(1-q)}\right)-f(1)\bigg|p_{\infty k}(q;qt)d_qt\\
&\leq&  [n+1]_q\int_0^1 q^{-k} \bigg|f\left(\frac{[n]_qt+\vartheta}{[n]_q+\varpi}\right)-f(1)\bigg| | p_{n k}(q;qt) -p_{\infty k}(q;qt) |d_qt \\
&+& \frac{1}{1-q}\int_0^1 q^{-k} \bigg| f\left(\frac{t+\vartheta(1-q)}{1+\varpi(1-q)}\right)-f(1)\bigg|p_{\infty k}(q;qt)d_qt+[n+1]_q\int_0^1 q^{-k} \bigg|f\left(\frac{[n]_qt+\vartheta}{[n]_q+\varpi}\right)-f(1)\bigg| p_{\infty k}(q;qt)d_qt  \\
 &\leq&  [n+1]_q\frac{q^{n-k}}{1-q}\int_0^1 q^{-k} \bigg|f\left(\frac{[n]_qt+\vartheta}{[n]_q+\varpi}\right)-f(1)\bigg| | p_{n k}(q;qt) +p_{\infty k}(q;qt) |d_qt \\
&+ &\frac{1}{1-q}\int_0^1 q^{-k} \bigg| f\left(\frac{t+\vartheta(1-q)}{1+\varpi(1-q)}\right)-f(1)\bigg|p_{\infty k}(q;qt)d_qt+[n+1]_q\int_0^1 q^{-k} \bigg|f\left(\frac{[n]_qt+\vartheta}{[n]_q+\varpi}\right)-f(1)\bigg| p_{\infty k}(q;qt)d_qt  \\
&\leq&  \omega(f,q^n) \left[2\frac{q^{n-k}}{1-q} \left(1+
q^{k+1-n} + \frac{q^{-n}(\varpi-\vartheta)}{[n]_q+\varpi}\right)+
\left(1+ q^{k+1-n} +
\frac{q^{-n}(\varpi-\vartheta)}{1+\varpi(1-q)}\right) + \left(1+
q^{k+1-n} +
\frac{q^{-n}(\varpi-\vartheta)}{[n]_q+\varpi}\right)\right].
\end{eqnarray*}
 To estimate $I_1, I_2$ and $I_3$, we have
\begin{eqnarray*}
I_1&\leq &
\frac{\omega(f,q^n)}{1-q}\left(8+\frac{3(\varpi-\vartheta)}{q^n([n]_q+\varpi)}\
                +\frac{(\varpi-\vartheta)}{q^{n}(1+\varpi(1-q))} \right) \sum_{k=0}^n p_{n k}(q;x)\\
                &=& \frac{\omega(f,q^n)}{1-q}\left(8+\frac{3(\varpi-\vartheta)}{q^n([n]_q+\varpi)}+\frac{(\varpi-\vartheta)}{q^{n}(1+\varpi(1-q))} \right);
\end{eqnarray*}
\begin{eqnarray*}
I_3&=&\sum_{k=n+1}^{\infty}| A_{\infty k
}^{\vartheta,\varpi}(f-f(1))|
p_{\infty k}(q;x)\\
&\leq& \omega(f,q^n) \sum_{k=n+1}^{\infty} \left(1+ q^{k+1-n}+ \frac{q^{-n}(\varpi-\vartheta)}{1+\varpi(1-q)}\right)p_{\infty k}(q;x)\\
 &\leq& \omega(f,q^n)\left(2+\frac{q^{-n}(\varpi-\vartheta)}{1+\varpi(1-q)}\right);\\
I_2&=&\sum_{k=0}^{n} A_{\infty k}^{\vartheta,\varpi} \left(f-f(1)\right) |p_{n k}(q;x)- p_{\infty k}(q;x)|\\
&\leq& \sum_{k=0}^{n} \left[ \omega(f,q^n)\left(1 +
q^{k+1-n}+\frac{q^{-n}(\varpi-\vartheta)}{1+\varpi(1-q)}\right)\right]\left[\frac{q^{n-k}}{1-q}|
p_{n k}(q;x)+p_{ \infty k}(q;x)|\right]\\
&\leq&
\frac{2\omega(f,q^n)}{1-q}\left(2+\frac{q^{-n}(\varpi-\vartheta)}{1+\varpi(1-q)}\right).
\end{eqnarray*}
Combining the estimates $I_1 - I_3$, we conclude that $\displaystyle \| D_{n,q }^{\vartheta,\varpi}(f)-D_{\infty,q}^{\vartheta,\varpi}(f)\| \leq C_q^{\vartheta,\varpi}\omega(f,q^n).$\\
 This complete the proof of Theorem \ref{t3}.
\begin{lemma}[\cite{zeng2010note}]\label{l3}
Let $L$ be a positive linear operator on $C\left([0,1]\right)$
which reproduces constant functions. If $L(t,x)>x \text{  }\forall
\text{  } x\in (0,1)$, then $L(f)=f$ if and only if $f$ is
constant.
\end{lemma}
\begin{remark}
Since $\displaystyle D_{\infty,q}^{\vartheta,\varpi}(t;x)=
\frac{(1+q(x-1))+\vartheta(1-q)}{1+\varpi(1-q)}>x $ for $0<q<1$
consequence of Lemma \ref{l3}, we have  the following:
\end{remark}
\begin{theorem}\label{t6}
Let $0 < q < 1$ be fixed and let $f \in C\left([0, 1]\right)$.
Then $D_{\infty,q}^{\vartheta,\varpi} (f; x) = f(x)$ for all $x\in
[0,1]$ if and only if $f$ is constant.
\end{theorem}
\begin{theorem}\label{t5}
For any $f\in C\left([0,1]\right)$, $\{D_{\infty, q}^{\vartheta,
\varpi}(f)\}$ converges to $f$ uniformly on $[0,1]$ as $q
\rightarrow 1-$.
\end{theorem}
\textbf{Proof:} We know that the operators $ D_{\infty,
q}^{\vartheta, \varpi}$ is positive linear operator on
$C\left([0,1]\right)$ and reproduce constant functions. Also,
$\displaystyle D_{\infty,q}^{\vartheta,\varpi}(t;x)\displaystyle
\rightarrow x \text{ uniformly on $[0,1]$ as } q\rightarrow 1-$
and $\displaystyle
D_{\infty, q}^{\vartheta, \varpi}(t^2;x)\displaystyle \rightarrow x^2 \text{ uniformly on $[0,1]$ as } q\rightarrow 1-.$\\
Thus, Theorem \ref{t5} follows from Korovkin Theorem.


\section*{References:}
\biboptions{numbers,sort&compress}

\end{document}